\documentclass{amsart}
\usepackage{latexsym}
\usepackage{amsfonts}
\usepackage{amssymb}
\usepackage{graphicx}
\usepackage[all]{xy}

\newtheorem{theorem}{Theorem}[section]
\newtheorem{corollary}{Corollary}[section]
\newtheorem{lemma}{Lemma}[section]
\newtheorem{sublemma}{Sublemma}[lemma]
\newtheorem{proposition}{Proposition}[section]
\newtheorem{fact}{Fact}[section]
\theoremstyle{remark}
\newtheorem{remark}{Remark}[section]

\numberwithin{equation}{section}

\newcommand{\g}{\mathfrak{g}}
\newcommand{\ak}{\mathfrak{k}}
\newcommand{\ba}{\mathfrak{a}}
\newcommand{\ab}{\mathfrak{b}}
\newcommand{\ah}{\mathfrak{h}}
\newcommand{\ap}{\mathfrak{p}}

\newcommand{\am}{\mathfrak{m}}
\newcommand{\at}{\mathfrak{t}}
\newcommand{\ZZ}{\mathbb{Z}}
\newcommand{\RR}{\mathbb{R}}
\newcommand{\Zn}{\mathbb{Z}/n\mathbb{Z}}
\newcommand{\Aut}{\mathrm{Aut}}
\newcommand{\Ad}{\mathrm{Ad}}
\newcommand{\Map}{\mathrm{Map}}
\newcommand{\Inn}{\mathrm{Inn}}

\begin{document}

\title{Nonabelian cohomology with coefficients in Lie
groups}

\author{Jinpeng An}
\address{School of mathematical science, Peking University,
 Beijing, 100871, P. R. China }
\email{anjinpeng@gmail.com}

\author{Zhengdong Wang}
\address{School of mathematical science, Peking University,
 Beijing, 100871, P. R. China}
\email{zdwang@pku.edu.cn}

\thanks{This work is supported by the 973 Project
Foundation of China (\#TG1999075102).}

\keywords{nonabelian cohomology, Lie group, twisted conjugate
action.}

\subjclass[2000]{20J06; 22E15; 57S15; 57S20.}

\begin{abstract}
In this paper we prove some properties of the nonabelian
cohomology $H^1(A,G)$ of a group $A$ with coefficients in a
connected Lie group $G$. When $A$ is finite, we show that for
every $A$-submodule $K$ of $G$ which is a maximal compact subgroup
of $G$, the canonical map $H^1(A,K)\rightarrow H^1(A,G)$ is
bijective. In this case we also show that $H^1(A,G)$ is always
finite. When $A=\ZZ$ and $G$ is compact, we show that for every
maximal torus $T$ of the identity component $G_0^\ZZ$ of the group
of invariants $G^\ZZ$, $H^1(\ZZ,T)\rightarrow H^1(\ZZ,G)$ is
surjective if and only if the $\ZZ$-action on $G$ is
$1$-semisimple, which is also equivalent to that all fibers of
$H^1(\ZZ,T)\rightarrow H^1(\ZZ,G)$ are finite. When $A=\Zn$, we
show that $H^1(\Zn,T)\rightarrow H^1(\Zn,G)$ is always surjective,
where $T$ is a maximal compact torus of the identity component
$G_0^{\Zn}$ of $G^{\Zn}$. When $A$ is cyclic, we also interpret
some properties of $H^1(A,G)$ in terms of twisted conjugate
actions of $G$.
\end{abstract}

\maketitle

%%%%%%%%%%%%%%%%%%%%%%%%%%%%%%%%%%%%%%%%%%%%%%%%%%%%%%%%%%% Section 1

\section{Introduction}

Let $A$ be a group, and let $G$ be a (nonabelian) $A$-module. A
map $\alpha:A\rightarrow G$ $(a\mapsto g_a)$ is called a cocycle
if $g_{ab}=g_aa(g_b)$. The set of all cocycles is denoted by
$Z^1(A,G)$. Two cocycles $\alpha_1,\alpha_2$ are cohomologous if
there exists $g\in G$ such that
$\alpha_2(a)=g^{-1}\alpha_1(a)a(g)$ for every $a\in A$. This is an
equivalence relation in $Z^1(A,G)$. The nonabelian cohomology
$H^1(A,G)$ of $A$ with coefficients in $G$ is defined to be the
set of all equivalence classes in $Z^1(A,G)$.

Most studies of this kind of cohomology concentrate on the case
that $G$ is an algebraic group, which play an important role in
the study of algebraic groups and number theory (see, for example,
\cite{NSW} and \cite{PR}). In this paper we consider the case that
$G$ is a connected Lie group, and prove some properties of
$H^1(A,G)$. Some of these properties provide a way to compute
$H^1(A,G)$. When $A$ is a cyclic group, we also relate $H^1(A,G)$
with the so-called twisted conjugate actions of $G$, and interpret
some properties of $H^1(A,G)$ as properties concerning structures
of orbits of twisted conjugate actions.

Throughout this paper, we make the convention that when we say $G$
is a Lie group with a nonabelian $A$-module structure, we always
assume that $A$ acts on $G$ by automorphisms of Lie group, that
is, $A$ acts smoothly on $G$; and when we say that $f:H\rightarrow
G$ is an $A$-homomorphism of such two $A$-modules, we always
assume that $f$ is also a homomorphism of Lie groups, that is, $f$
is smooth.

Two cases regarding the type of the group $A$ will be considered,
that is, $A$ is a finite group and $A$ is the infinite cyclic
group $\ZZ$. For the first case, if $A$ is a finite cyclic group
$\Zn$, further properties of $H^1(\Zn,G)$ will be proved.

For the case that $A$ is finite, we have

\begin{theorem}\label{A-bi}
Let $A$ be a finite group, and let $G$ be a connected Lie group
with an
$A$-module structure. Then\\
(i) There exists a maximal compact subgroup $K$ of $G$ which is an
$A$-submodule of $G$;\\
(ii) For each maximal compact subgroup $K$ of $G$ which is also an
$A$-submodule of $G$, the canonical map $H^1(A,K)\rightarrow
H^1(A,G)$ is bijective.
\end{theorem}

Theorem \ref{A-bi} reduces the calculation of $H^1(A,G)$ to the
calculation of the cohomology of $A$ with coefficients in a
compact Lie group. In the case that $G$ is the complexification of
a compact Lie group $K$ and $\ZZ/2\ZZ$ acts on $G$ by complex
conjugation, the bijectivity of the map
$H^1(\ZZ/2\ZZ,K)\rightarrow H^1(\ZZ/2\ZZ,G)$ has been proved in
Serre \cite{Se2}, Chapter III, Section 4.5.

Let $\Map(A,G)$ denote the set of all maps from $A$ to $G$. If $A$
is finite, there is a natural identification $\Map(A,G)\cong
G^{|A|}$, where $|A|$ is the order of $A$. Hence $\Map(A,G)$
inherits a structure of smooth manifold from $G^{|A|}$. The set of
cocycles $Z^1(A,G)$ is a closed subset of $\Map(A,G)$. The
following assertion will be proved.

\begin{theorem}\label{A-cocycle}
Let $A$ be a finite group, and let $G$ be a connected Lie group
with an $A$-module structure. Then we have\\
(i) $Z^1(A,G)\subset\Map(A,G)$ has finitely many connected
components, each of which is a closed submanifold of $\Map(A,G)$;\\
(ii) Each cohomology class in $Z^1(A,G)$ is a connected component
of $Z^1(A,G)$.
\end{theorem}

Theorem \ref{A-cocycle} has the following obvious corollary.

\begin{corollary}\label{A-finite}
Let $A$ be a finite group, and let $G$ be a connected Lie group
with an $A$-module structure. Then $H^1(A,G)$ is finite.
\end{corollary}

Corollary \ref{A-finite} generalizes the finiteness theorem for
Galois cohomology with coefficients in algebraic groups defined
over $\RR$ (see \cite{Se2}, Chapter III, Section 4.3).

For the case that $A=\ZZ$, we restrict our attention to the case
that $G$ is compact. To state our theorem, we first introduce a
notion.

An automorphism $\sigma$ of a connected compact Lie group $G$ is
\emph{$1$-semisimple} if $\ker(1-d\sigma)=\ker((1-d\sigma)^2)$ in
the Lie algebra $\g$ of $G$. All automorphisms of a compact
semisimple Lie group are $1$-semisimple, and then $\sigma$ is
$1$-semisimple if and only if the restriction of $\sigma$ to the
connected component of the center of $G$ is $1$-semisimple. For a
connected compact Lie group $G$ with a $\ZZ$-module structure, the
$\ZZ$-action on $G$ is \emph{$1$-semisimple} if a generator of
$\ZZ$ acts $1$-semisimply on $G$.

For a Lie group $G$ with an $A$-module structure, we always denote
the group of invariants by $G^A$, and denote the identity
component of $G^A$ by $G_0^A$. The following theorem will be
proved.

\begin{theorem}\label{Z-sur}
Let $G$ be a connected compact Lie group with a $\ZZ$-module
structure, $T$ a maximal torus of $G_0^\ZZ$. Let $i:T\rightarrow
G$ be the canonical inclusion, and let $i_1:H^1(\ZZ,T)\rightarrow
H^1(\ZZ,G)$ be the canonical map induced by $i$. Then the
following statements are equivalent.\\
(i) The $\ZZ$-action on $G$ is $1$-semisimple;\\
(ii) The map $i_1$ is surjective;\\
(iii) The kernel of $i_1$ is finite;\\
(iv) All fibers of $i_1$ are finite.
\end{theorem}

We note here that when we say that a set is finite, we allow it to
be empty. We will use the word ``nonempty and finite" if we want
to emphasize that it is nonempty.

If $A$ is a finite cyclic group $\Zn$, we have further properties
about $H^1(A,G)$ besides Theorems \ref{A-bi} and \ref{A-cocycle}.

For a Lie group $H$, a subgroup $T$ of $H$ is a \emph{maximal
compact torus} of $H$ if $T$ is a compact torus and there is no
other compact torus $T'$ of $H$ such that $T\subsetneqq T'$. In
fact, any two maximal compact torus are conjugate (see Section 5).
The following assertion will be proved.

\begin{theorem}\label{Zn-sur}
Let $G$ be a connected Lie group with a $\Zn$-module structure,
$T$ a maximal compact torus of $G_0^{\Zn}$. Then the canonical map
$i_1:H^1(\Zn,T)\rightarrow H^1(\Zn,G)$ is surjective.
\end{theorem}

A theorem of T. A. Springer (see Serre \cite{Se2}, Chapter III,
Section 4.3, Lemma 6) says that if $G$ is a linear algebraic group
defined over a perfect field $k$ with a Cartan subgroup $C$, then
the canonical map $H^1(k,N_G(C))\rightarrow H^1(k,G)$ between sets
of Galois cohomology is surjective. Theorem \ref{Zn-sur} may be
viewed as a generalization of the theorem of Springer in the case
that $k=\RR$.

The set of cocycles $Z^1(\Zn,G)$ may be identified with a closed
subset $Z=\{g\in G|g\sigma(g)\cdots\sigma^{n-1}(g)=e\}$ of $G$,
where $\sigma$ is a generator of $\Zn$ (see Section 2). Under this
identification, Theorem \ref{A-cocycle} has the following variant
form.

\begin{theorem}\label{Zn-cocycle}
Let $G$ be a connected Lie group with a $\Zn$-module structure,
$\sigma$ a generator of $\Zn$. Let $Z=\{g\in
G|g\sigma(g)\cdots\sigma^{n-1}(g)=e\}$. Then we have\\
(i) $Z$ has finitely many connected components, each of which is a
closed submanifold of $G$;\\
(ii) Under the identification $Z^1(\Zn,G)\cong Z$, each cohomology
class is a connected component of $Z$.
\end{theorem}

When $A$ is cyclic, that is, $A=\ZZ$ or $\Zn$, the first
nonabelian cohomology $H^1(A,G)$ has closed relation with the
so-called twisted conjugate actions of $G$ on $G$. For a connected
Lie group $G$ with an automorphism $\sigma$, the \emph{twisted
conjugate action} $\tau:G\times G\rightarrow G$ of $G$ on $G$
associated with $\sigma$ is defined by
$\tau_g(h)=gh\sigma(g)^{-1}$. An orbit of the twisted conjugate
action is called a \emph{twisted orbit}. This kind of action was
introduced in \cite{AW1,AW2} where an equivariant embedding of a
symmetric space $G/K$ into $G$ (or a covering group of $G$) was
established, and the image of the embedding is just the twisted
orbit through the identity element of $G$. This embedding was then
applied to the theory of random matrix ensembles associated with
symmetric spaces in \cite{AWY1,AWY2}, transforming the integration
manifold of a random matrix ensemble from a symmetric space
modelled by a matrix group to a space of matrices.

The relation between twisted conjugate actions and nonabelian
cohomology may be described briefly as follows. An automorphism
$\sigma$ of a Lie group $G$ induces naturally a $\ZZ$-module
structure on $G$ by $(m,g)\mapsto\sigma^m(g)$, where $m\in\ZZ,
g\in G$. The set of cocycles $Z^1(\ZZ,G)$ may be identified
naturally with $G$, and two cocycles $z_1,z_2\in G$ are
cohomologous if and only if they lie in the same twisted orbit. So
$H^1(\ZZ,G)$ may be identified with the space of twisted orbits.
If moreover $\sigma$ has finite order dividing a positive integer
$n$, the associated $\ZZ$-module structure on $G$ reduces to a
$\Zn$-module structure on $G$. The set of cocycles $Z^1(\Zn,G)$
may be identified naturally with $Z=\{g\in
G|g\sigma(g)\cdots\sigma^{n-1}(g)=e\}$, which is invariant under
the twisted conjugate action. Two cocycles $z_1,z_2\in Z$ are
cohomologous if and only if they lie in the same twisted orbit.
Hence $H^1(\Zn,G)$ may be identified with the set of twisted
orbits which are contained in $Z$. More details about this
relation will be given in Section 2.

With this relation in hand, we will interpret some properties of
$H^1(A,G)$ mentioned above in the language of twisted conjugate
actions. Some of these interpretations are as follows.

\begin{theorem}\label{compact-twisted}
Let $G$ be a connected compact Lie group, $\sigma$ an automorphism
of $G$. Let $T$ be a maximal torus of $G_0^\sigma$. Then the
following statements are equivalent.\\
(i) $\sigma$ is $1$-semisimple;\\
(ii) $O\cap T$ is nonempty for every twisted orbit $O$ of $G$
associated
with $\sigma$;\\
(iii) $O_e\cap T$ is finite, where $O_e$ is the twisted orbit
through the identity element $e$;\\
(iv) $O\cap T$ is finite for every twisted orbit $O$ of $G$
associated with $\sigma$.
\end{theorem}

\begin{theorem}\label{noncompact-twisted}
Let $G$ be a connected Lie group, $\sigma$ an automorphism of $G$
of finite order. Let $Z=\{g\in
G|g\sigma(g)\cdots\sigma^{n-1}(g)=e\}$, where $n$ is a positive
integer which is divisible
by the order of $\sigma$. Then\\
(i) $Z$ has finitely many connected components, each of which is a
twisted orbit;\\
(ii) For every maximal compact torus $T$ of $G_0^\sigma$ and every
twisted orbit $O$ contained in $Z$, $O\cap T$ is nonempty and
finite;\\
(iii) There exists a maximal compact subgroup $K$ of $G$ which is
$\sigma$-invariant. For every such $K$ and every twisted orbit $O$
of $G$ contained in $Z$, $O\cap K$ is a twisted orbit of $K$.
\end{theorem}

Now we give a sketch of the contents of each of the following
sections. In Section 2 we will review some basic facts on
nonabelian cohomology. The relation between nonabelian cohomology
of cyclic groups and twisted conjugate actions will also be given.

Section 3 will concern properties of nonabelian cohomology of
finite groups with coefficients in Lie groups. Theorems \ref{A-bi}
and \ref{A-cocycle} will be proved. A key ingredient in the proof
of Theorem \ref{A-bi} is the conjugacy theorem for maximal compact
subgroups of Lie groups with finitely many connected components.
The proof of Theorem \ref{A-cocycle} makes use of a theorem of
Weil for smooth manifolds, which we call Weil's Lemma (Fact
\ref{F:Weil}).

In Section 4 we will prove Theorem \ref{Z-sur}. Note that Theorem
\ref{Z-sur} is equivalent to Theorem \ref{compact-twisted}, which
says, among other things, that if the automorphism $\sigma$ is
$1$-semisimple, then every twisted orbit intersects every maximal
torus of $G^\sigma_0$. This fact is a generalization of Cartan's
conjugacy theorem for compact Lie groups, which is the special
case that $\sigma$ is the identity in Theorem
\ref{compact-twisted}. Besides using cohomology exact sequences
and twisting processes, an important ingredient in the proof of
Theorem \ref{Z-sur} is the Lefschetz Fixed Point Theorem. This is
reflected in the proof of Lemma \ref{L:ss-sur}. Our proof of Lemma
\ref{L:ss-sur} follows the idea of Weil's topological proof of
Cartan's conjugacy theorem \cite{We}. We also show that under the
condition of $1$-semisimplicity, each connected component of
$H^0(\ZZ,G/T)$ is a closed homogeneous submanifold of $G/T$, and
the restriction of the coboundary operator
$\delta:H^0(\ZZ,G/T)\rightarrow H^1(\ZZ,T)$ to each connected
component of $H^0(\ZZ,G/T)$ is constant.

Theorems \ref{Zn-sur} and \ref{Zn-cocycle} will be proved in
Section 5. The proof of Theorem \ref{Zn-sur} will make use of
Theorems \ref{A-bi} and \ref{Z-sur}. A corollary of Theorem
\ref{Zn-sur} says that if $G$ is simply connected and solvable,
then $H^1(\Zn,G)$ is trivial, due to the fact that $G$ admits no
nontrivial compact torus. Theorem \ref{Zn-cocycle} is in fact
equivalent to the particular case of Theorem \ref{A-cocycle} when
$A$ is finite cyclic. We will provide another proof of Theorem
\ref{Zn-cocycle} in Section 5, without the use of Weil's Lemma. We
also have the result that each connected component of
$H^0(\Zn,G/T)$ is a closed homogeneous submanifold of $G/T$, and
the restriction of the coboundary operator
$\delta:H^0(\Zn,G/T)\rightarrow H^1(\Zn,T)$ to each connected
component of $H^0(\Zn,G/T)$ is constant.

Section 6 will be devoted to properties of twisted conjugate
actions of Lie groups. We will first present some elementary
properties of such actions. Then we interpret some results proved
in Sections 3--5 in the language of twisted conjugate actions.
Theorems \ref{compact-twisted} and \ref{noncompact-twisted} will
be proved. In fact, Theorem \ref{compact-twisted} is equivalent to
Theorem \ref{Z-sur}, item (i) of Theorem \ref{noncompact-twisted}
is equivalent to Theorem \ref{Zn-cocycle}, item (ii) of Theorem
\ref{noncompact-twisted} can be implied from Corollary
\ref{A-finite} and Theorem \ref{Zn-sur}, item (iii) of Theorem
\ref{noncompact-twisted} can be implied from Theorem \ref{A-bi}.

The authors would like to thank the two referees, for suggestions
on improving the presentation of the content (especially the
comment of considering cohomology of general finite groups) from
one of them, and for suggestions on improving the English usage
from the other. The first author is in debt to Professor Jiu-Kang
Yu for many valuable conversations which engaged the author in the
subject of nonabelian cohomology. He also would like to thank
Professors F. Fang, Z. Hajto, and K.-H. Neeb for valuable
conversations or kind help.

%%%%%%%%%%%%%%%%%%%%%%%%%%%%%%%%%%%%%%%%%%%%%%%%%%%%%%%  Section 2

\section{Preliminaries on nonabelian cohomology}

In this section we list some facts of nonabelian cohomology which
will be used later, and explain the relation between nonabelian
cohomology and twisted conjugate actions. Most of these facts may
be found in Serre \cite{Se1,Se2}.

Let $A$ be a group, $X$ an $A$-set. By definition, the zeroth
nonabelian cohomology $H^0(A,X)$ of $A$ with coefficients in $X$
is the set $X^A=\{x\in X|a(x)=x, \forall a\in A\}$. Now let $G$ be
a group (which is always a Lie group in subsequent sections) on
which $A$ acts by automorphisms, that is, $G$ is a nonabelian
$A$-module, then $H^0(A,G)=G^A$ is a subgroup of $G$. A map
$\alpha:A\rightarrow G$ $(a\mapsto g_a)$ is called a cocycle if
$g_{ab}=g_aa(g_b)$ (This forces that $\alpha$ maps the identity
element of $A$ to the identity element of $G$). The set of all
cocycles is denoted by $Z^1(A,G)$. Two cocycles
$\alpha_1,\alpha_2\in Z^1(A,G)$ are cohomologous if there exists
$g\in G$ such that $\alpha_2(a)=g^{-1}\alpha_1(a)a(g)$ for every
$a\in A$. This is an equivalence relation in $Z^1(A,G)$, and the
first nonabelian cohomology $H^1(A,G)$ of $A$ with coefficients in
$G$ is defined to be the pointed set of all equivalence classes in
$Z^1(A,G)$ with neutral element the class of the unit cocycle. It
is obvious that $H^1(A,G\times H)\cong H^1(A,G)\times H^1(A,H)$.
Note that for nonabelian modules, only the zeroth and first
cohomology are commonly used. For the definition of higher
cohomology, see the references in \cite{Se1,Se2}.

Here is an equivalent description of $H^1(A,G)$. Let
$\alpha:A\rightarrow G$ $(a\mapsto g_a)$ be a cocycle. Construct a
map $\widetilde{\alpha}:A\rightarrow G\rtimes A$ by
$\widetilde{\alpha}(a)=(g_a,a)$. Then $\widetilde{\alpha}$ is
homomorphism of groups. Conversely, if
$\widetilde{\alpha}:A\rightarrow G\rtimes A$ is a homomorphism
splitting the exact sequence $0\rightarrow G\rightarrow G\rtimes
A\rightarrow A\rightarrow 0$, then $p_1\circ \widetilde{\alpha}$
is a cocycle, where $p_1:G\rtimes A\rightarrow G$ is the
projection to the first factor. So $Z^1(A,G)$ may be identified
with the set of homomorphisms from $A$ to $G\rtimes A$ splitting
the exact sequence $0\rightarrow G\rightarrow G\rtimes
A\rightarrow A\rightarrow 0$. Under this identification, it is
easy to show that two homomorphisms
$\widetilde{\alpha}_1,\widetilde{\alpha}_2:A\rightarrow G\rtimes
A$ are cohomologous if and only if there exists $g\in G$ such that
$\widetilde{\alpha}_2=g^{-1}\widetilde{\alpha}_1g$.

Let $f:H\rightarrow G$ be a homomorphism of nonabelian
$A$-modules. Then there are canonical induced maps of pointed sets
$f_i:H^i(A,H)\rightarrow H^i(A,G), i=0,1$. $f_0$ is also a
homomorphism of groups. If $H$ is an $A$-submodule of $G$, then
$G/H$ is an $A$-set, and we have an exact sequence of pointed sets
$$
0\rightarrow H^0(A,H)\stackrel{f_0}\rightarrow H^0(A,G)\rightarrow
H^0(A,G/H)\stackrel{\delta}\rightarrow
H^1(A,H)\stackrel{f_1}\rightarrow H^1(A,G),
$$
where the coboundary operator $\delta$ maps an element $gH\in
H^0(A,G/H)$ to the class of the cocycle $g_a=g^{-1}a(g)$, and the
neutral element of $H^0(A,G/H)$ is chosen to be $H\in G/H$.

Let $G$ be a nonabelian $A$-module, and let $\alpha\in Z^1(A,G)$.
One can define another $A$-module $G_\alpha$ by twisting $G$ using
the cocycle $\alpha$ as follows. The underlying group of
$G_\alpha$ is $G$, and $A$ acts on $G_\alpha$ by
$t_a(g)=\alpha(a)a(g)\alpha(a)^{-1}$. Then the map
$Z^1(A,G_\alpha)\rightarrow Z^1(A,G)$ defined by
$\beta\mapsto\beta\cdot\alpha$ is a bijection, and induces a
bijection $t_\alpha:H^1(A,G_\alpha)\rightarrow H^1(A,G)$, under
which the neutral element of $H^1(A,G_\alpha)$ is mapped to the
cohomology class $[\alpha]$ of $\alpha$. If $f:H\rightarrow G$ is
an $A$-homomorphism of nonabelian $A$-modules, and if $\alpha\in
Z^1(A,H)$, then $f_\alpha=f:H_\alpha\rightarrow G_{f\circ\alpha}$
is also an $A$-homomorphism, and we have the following commutative
diagram
$$
\xymatrix{H^1(A,H_\alpha)\ar[d]^{t_\alpha}\ar[r]^{f_{\alpha1}}
& H^1(A,G_{f\circ\alpha})\ar[d]^{t_{f\circ\alpha}}\\
H^1(A,H) \ar[r]^{f_1} & H^1(A,G).}
$$
So $t_\alpha$ is a bijection between $\ker(f_{\alpha1})$ and the
fiber $f_1^{-1}(f_1([\alpha]))$. This means that the twisting
process allows one to transform each non-empty fiber of $f_1$ to a
kernel.

Now we describe the relation between nonabelian cohomology and
twisted conjugate actions. Let $G$ be a group with an automorphism
$\sigma$. the twisted conjugate action $\tau$ of $G$ on itself
associated with $\sigma$ is defined by
$\tau_g(h)=gh\sigma(g)^{-1}$. An orbit of the twisted conjugate
action is called a twisted orbit. The automorphism $\sigma$ of $G$
induces in a natural way a $\ZZ$-action $\ZZ\times G\rightarrow G$
on $G$, defined by $(m,g)\mapsto\sigma^m(g)$. This makes $G$ a
nonabelian $\ZZ$-module. A cocycle in $Z^1(\ZZ,G)$ is determined
by its value on $1\in\ZZ$. In fact, for a cocycle $\ZZ\rightarrow
G$ $(m\mapsto g_m)$, we have
$g_m=g_1\sigma(g_1)\cdots\sigma^{m-1}(g_1)$,
$g_{-m}=\sigma^{-1}(g_1^{-1})\cdots\sigma^{-m}(g_1^{-1}), m\geq
1$. Conversely, each $g_1\in G$ determines a cocycle in this way.
So we may identify $Z^1(\ZZ,G)$ with $G$. Under this
identification, two cocycles $z_1,z_2\in G$ are cohomologous if
and only if there exists $g\in G$ such that
$z_2=g^{-1}z_1\sigma(g)$, that is, $z_1$ and $z_2$ lie in the same
twisted orbit. So $H^1(\ZZ,G)$ may be identified with the space of
twisted orbits.

Let $H$ be a subgroup of $G$ which is $\sigma$-invariant. Then $H$
is a $\ZZ$-submodule of $G$. With the identification
$Z^1(\ZZ,G)\cong G$ in mind, it is clear that to say the canonical
map $H^1(\ZZ,H)\rightarrow H^1(\ZZ,G)$ is surjective is equivalent
to say that the intersection of every twisted orbit of $G$ with
$H$ is nonempty, and to say $H^1(\ZZ,H)\rightarrow H^1(\ZZ,G)$ is
bijective is equivalent to say that the intersection of every
twisted orbit of $G$ with $H$ is just one twisted orbit of $H$. In
particular, if $H\subset G^\sigma=\{g\in G|\sigma(g)=g\}$ and $H$
is abelian, then $H^1(\ZZ,H)=Z^1(\ZZ,H)$, hence to say the kernel
of $H^1(\ZZ,H)\rightarrow H^1(\ZZ,G)$ is finite is equivalent to
say that the intersection of the twisted orbit through the
identity element $e$ of $G$ with $H$ is a finite set, and to say
all fibers of $H^1(\ZZ,H)\rightarrow H^1(\ZZ,G)$ are finite is
equivalent to say that the intersection of every twisted orbit
with $H$ is a finite set.

If $\sigma$ is of finite order, the associated $\ZZ$-module
structure on $G$ reduces to a $\Zn$-module structure on $G$ for
each positive integer $n$ which is divisible by the order of
$\sigma$. Let $\Zn\rightarrow G$ $(m\mapsto g_m)$ be a cocycle in
$Z^1(\Zn,G)$. Then we have
$g_m=g_1\sigma(g_1)\cdots\sigma^{m-1}(g_1), 1\leq m\leq n-1$, and
$g_1\sigma(g_1)\cdots\sigma^{n-1}(g_1)=e$. Conversely, each
$g_1\in G$ satisfying $g_1\sigma(g_1)\cdots\sigma^{n-1}(g_1)=e$
determines a cocycle in this way. So we may identify $Z^1(\Zn,G)$
with the subset $Z=\{g\in G|g\sigma(g)\cdots\sigma^{n-1}(g)=e\}$
of $G$, which is invariant under the twisted conjugate action.
Under this identification, two cocycles $z_1,z_2\in Z$ are
cohomologous if and only if they lie in the same twisted orbit.
Hence $H^1(\Zn,G)$ may be identified with the set of twisted
orbits which are contained in $Z$.

If $H$ is a closed subgroup of $G$ which is $\sigma$-invariant,
then $H$ is a $\Zn$-submodule of $G$. Similar to the case of
$\ZZ$-modules, to say the canonical map $H^1(\Zn,H)\rightarrow
H^1(\Zn,G)$ is surjective is equivalent to say that the
intersection of every twisted orbit contained in $Z$ with $H$ is
nonempty, and so on.

%%%%%%%%%%%%%%%%%%%%%%%%%%%%%%%%%%%%%%%%%%%%%%%%%%%%%%%  Section 3

\section{Nonabelian cohomology of finite groups with coefficients in Lie
groups}

In this section we consider nonabelian cohomology of general
finite groups with coefficients in connected Lie groups, and prove
Theorems \ref{A-bi} and \ref{A-cocycle}. The first theorem reads
as follows.

\begin{theorem}\label{T:A-bi}
Let $A$ be a finite group, $G$ a connected Lie group with an
$A$-module structure. Then\\
(i) There exists a maximal compact subgroup $K$ of $G$ which is an
$A$-submodule of $G$;\\
(ii) For each $K$ satisfying the conditions in (i), the canonical
map $H^1(A,K)\rightarrow H^1(A,G)$ is bijective.
\end{theorem}

To prove Theorem \ref{T:A-bi}, we need the following fact, which
will also be used later in the proofs of other assertions.

\begin{fact}[\cite{Ho}, Chapter XV, Theorem 3.1]\label{F:max-cpt}
Let $H$ be a Lie group with finitely many connected components,
$K$ a maximal compact
subgroup of $H$. Let $\ah$ and $\ak$ be the Lie algebras of $H$
and $K$, respectively. Then we have\\
(i) $K\cap H_0=K_0$, and $K_0$ is a maximal compact
subgroup of $H_0$;\\
(ii) For each compact subgroup $K'$ of $H$, there exists $h\in
H_0$ such that $hK'h^{-1}\subset K$;\\
(iii) there exist linear subspaces $\am_1,\cdots,\am_r$ of $\ah$
with $\ah=\ak\oplus\am_1\oplus\cdots\oplus\am_r$ such that
$\Ad(k)(\am_i)=\am_i, \forall k\in K, i\in\{1,\cdots,r\}$, and
such that the map
$\varphi:K\times\am_1\times\cdots\times\am_r\rightarrow H$ defined
by $\varphi(k,X_1,\cdots,X_r)=ke^{X_1}\cdots e^{X_r}$ is a
diffeomorphism.
\end{fact}

\begin{proof} [Proof of Theorem \ref{T:A-bi}]
We construct the semidirect product $G\rtimes A$, which is a Lie
group with finitely many connected components. We identify the
identity component of $G\rtimes A$ with $G$. Then
$a(g)=(e,a)g(e,a)^{-1}$ for all $a\in A$ and $g\in G$. Since the
subgroup $\{e\}\rtimes A$ of $G\rtimes A$ is compact, there is a
maximal compact subgroup $K'$ of $G\rtimes A$ such that
$\{e\}\rtimes A\subset K'$. Let $K=K'_0\subset G$. By Fact
\ref{F:max-cpt}, $K$ is a maximal compact subgroup of $G$. Since
$a(K)=(e,a)K(e,a)^{-1}=K$ for all $a\in A$, $K$ is an
$A$-submodule of $G$. This proves (i).

Now we prove (ii). Let $K$ be a maximal compact subgroup $K$ of
$G$ which is also a $\Zn$-submodule of $G$. Then $K\rtimes A$ is a
maximal compact subgroup of $G\rtimes A$. Note that $Z^1(A,G)$ may
be identified naturally with the set of homomorphisms
$A\rightarrow G\rtimes A$ splitting the exact sequence
$$0\rightarrow G\rightarrow G\rtimes
A\rightarrow A\rightarrow 0$$ by identifying a cocycle
$\alpha:A\rightarrow G$ with $\widetilde{\alpha}:A\rightarrow
G\rtimes A$ $(a\mapsto(\alpha(a),a))$, and $\alpha$ is in
$Z^1(A,K)$ if and only if $\widetilde{\alpha}$ assumes values in
$K\rtimes A$ (see Section 2). Hence to prove $H^1(A,K)\rightarrow
H^1(A,G)$ is surjective, it is sufficient to prove that for every
homomorphism $\widetilde{\alpha}:A\rightarrow G\rtimes A$
splitting the above exact sequence, there exists $g\in G$ such
that $g^{-1}\widetilde{\alpha}g$ assumes values in $K\rtimes A$.
But for such an $\widetilde{\alpha}$, $\widetilde{\alpha}(A)$ is a
finite subgroup of $G\rtimes A$. By Fact \ref{F:max-cpt}, there
exists $g\in G$ such that $g^{-1}\widetilde{\alpha}(A)g\subset
K\rtimes A$. This proves that $H^1(A,K)\rightarrow H^1(A,G)$ is
surjective.

Let $\g$, $\ak$ be the Lie algebras of $G$ and $K$. By Fact
\ref{F:max-cpt}, there exist linear subspaces $\am_1,\cdots,\am_r$
of $\g$ with $\g=\ak\oplus\am_1\oplus\cdots\oplus\am_r$ such that
$\Ad(k)(\am_i)=\am_i, \forall k\in K\rtimes A, 1\leq i\leq r$, and
such that the map
$\varphi:K\times\am_1\times\cdots\times\am_r\rightarrow G$ defined
by $\varphi(k,X_1,\cdots,X_r)=ke^{X_1}\cdots e^{X_r}$ is a
diffeomorphism. To prove the injectivity of $H^1(A,K)\rightarrow
H^1(A,G)$, let $\widetilde{\alpha}_1,
\widetilde{\alpha}_2:A\rightarrow K\rtimes A$ be two homomorphisms
splitting the above exact sequence. Suppose there exists $g\in G$
such that $\widetilde{\alpha}_2=g^{-1}\widetilde{\alpha}_1g$. It
is sufficient to show that there exists $k\in K$ such that
$\widetilde{\alpha}_2=k^{-1}\widetilde{\alpha}_1k$. For every
$a\in A$, rewrite the equality
$\widetilde{\alpha}_2(a)=g^{-1}\widetilde{\alpha}_1(a)g$ as
$\widetilde{\alpha}_2(a)^{-1}g\widetilde{\alpha}_2(a)=\widetilde{\alpha}_2(a)^{-1}\widetilde{\alpha}_1(a)g$,
and let $g=ke^{X_1}\cdots e^{X_r}, k\in K, X_i\in\am_i$. Then we
have
$$\widetilde{\alpha}_2(a)^{-1}k\widetilde{\alpha}_2(a)e^{\Ad(\widetilde{\alpha}_2(a)^{-1})X_1}\cdots
e^{\Ad(\widetilde{\alpha}_2(a)^{-1})X_r}
=\widetilde{\alpha}_2(a)^{-1}\widetilde{\alpha}_1(a)ke^{X_1}\cdots
e^{X_r},$$ that is,
\begin{align*}
&\varphi(\widetilde{\alpha}_2(a)^{-1}k\widetilde{\alpha}_2(a),\Ad(\widetilde{\alpha}_2(a)^{-1})X_1,\cdots,
\Ad(\widetilde{\alpha}_2(a)^{-1})X_r)\\
=&\varphi(\widetilde{\alpha}_2(a)^{-1}\widetilde{\alpha}_1(a)k,X_1,\cdots,X_r).
\end{align*}
Since $\varphi$ is a diffeomorphism, we have
$\widetilde{\alpha}_2(a)^{-1}k\widetilde{\alpha}_2(a)=\widetilde{\alpha}_2(a)^{-1}\widetilde{\alpha}_1(a)k,
\forall a\in A$, that is,
$\widetilde{\alpha}_2=k^{-1}\widetilde{\alpha}_1k$. This proves
the injectivity, and then finishes the proof of the theorem.
\end{proof}

For two sets $X$ and $Y$, denote the set of all maps from $X$ to
$Y$ by $\Map(X,Y)$. Let $A$ be a finite group, and let $G$ be a
connected Lie group with an $A$-module structure. Then there is a
natural identification $\Map(A,G)\cong\prod_{a\in A}G\cong
G^{|A|}$, where $|A|$ is the order of $A$. Hence $\Map(A,G)$
inherits a structure of smooth manifold from $G^{|A|}$. The set of
cocycles $Z^1(A,G)$ is a closed subset of $\Map(A,G)$.

\begin{theorem}\label{T:A-cocycle}
Let $A$ be a finite group, $G$ a connected Lie group
with an $A$-module structure. Then we have\\
(i) $Z^1(A,G)\subset\Map(A,G)$ has finitely many connected
components, each of which is a closed submanifold of $\Map(A,G)$;\\
(ii) Each cohomology class in $Z^1(A,G)$ is a connected component
of $Z^1(A,G)$, that is, $H^1(A,G)$ coincides with the set of all
connected components of $Z^1(A,G)$.
\end{theorem}

The following fact is needed in the proof of Theorem
\ref{T:A-cocycle}. For a proof of it, see \cite{Se3}, Part II,
Chapter III, Section 11.

\begin{fact}[Weil's Lemma]\label{F:Weil}
Let $M, N$ be smooth manifolds, $\varphi:M\rightarrow N$ a smooth
map. Let $x\in N$. Suppose that for every $y\in \varphi^{-1}(x)$,
there exist a smooth manifold $L$, a smooth map $\psi:L\rightarrow
M$, and a point $z\in
L$ such that\\
(\textit{i}) $\psi(z)=y$;\\
(\textit{ii}) $\varphi(\psi(w))=x$ for all $w\in L$;\\
(\textit{iii}) The sequence
$T_zL\stackrel{(d\psi)_z}\longrightarrow
T_yM\stackrel{(d\varphi)_y}\longrightarrow T_xN$ is exact.\\
Then every connected component of $\varphi^{-1}(x)$ is a closed
submanifold of $M$. Moreover, there exist a neighborhood $U$ of
$y$ in $M$ and a neighborhood $V$ of $z$ in $L$ such that
$\psi(V)=U\cap \varphi^{-1}(x)$.
\end{fact}

\begin{proof} [Proof of Theorem \ref{T:A-cocycle}]
Let $\g$ be the Lie algebra of $G$. Then for $\alpha\in\Map(A,G)$,
the tangent space $T_\alpha\Map(A,G)$ of $\Map(A,G)$ at $\alpha$
can be naturally identified with $\Map(A,\g)$, by the
correspondence that $\lambda\in\Map(A,\g)$ is mapped to
$\frac{d}{dt}\big|_{t=0}e^{t\lambda}\alpha\in T_\alpha\Map(A,G)$.
Equivalently, for a smooth curve $\gamma:\RR\rightarrow\Map(A,G)$
with $\gamma(0)=\alpha$, the tangent vector
$\frac{d}{dt}\big|_{t=0}\gamma(t)$ corresponds to the element of
$\Map(A,\g)$ defined by
$a\mapsto\frac{d}{dt}\big|_{t=0}\gamma(t)(a)\alpha(a)^{-1}$.
Similarly, tangent spaces of $\Map(A^2,G)$ may be identified with
$\Map(A^2,\g)$.

We first prove that every connected component of
$Z^1(A,G)\subset\Map(A,G)$ is a closed submanifold. Define a map
$\varphi:\Map(A,G)\rightarrow\Map(A^2,G)$ by
$\varphi(\alpha)(a,b)=\alpha(a)a(\alpha(b))\alpha(ab)^{-1}$.
Denote the constant map $(a,b)\mapsto e$ in $\Map(A^2,G)$ also by
$e$. Then $\varphi^{-1}(e)=Z^1(A,G)$. Now for $\alpha$ $(a\mapsto
g_a)\in Z^1(A,G)$, we compute $\ker((d\varphi)_\alpha)$. For
$\lambda$ $(a\mapsto X_a)\in\Map(A,\g)\cong T_\alpha\Map(A,G)$ and
$a,b\in A$, we have
\begin{align*}
&(d\varphi)_\alpha(\lambda)(a,b)\\
=&(d\varphi)_\alpha\left(\frac{d}{dt}\Big|_{t=0}e^{t\lambda}\alpha\right)(a,b)\\
=&\frac{d}{dt}\Big|_{t=0}\varphi(e^{t\lambda}\alpha)(a,b)\\
=&\frac{d}{dt}\Big|_{t=0}(e^{tX_a}g_a)a(e^{tX_b}g_b)(e^{tX_{ab}}g_{ab})^{-1}\\
=&\frac{d}{dt}\Big|_{t=0}e^{tX_a}(g_ae^{tda(X_b)}a(g_b)g_{ab}^{-1})e^{-tX_{ab}}\\
=&\frac{d}{dt}\Big|_{t=0}e^{tX_a}e^{t\Ad(g_a)da(X_b)}e^{-tX_{ab}}\\
=&X_a+\Ad(g_a)da(X_b)-X_{ab}.
\end{align*}
So $\lambda$ $(a\mapsto X_a)\in\ker((d\varphi)_\alpha)$ if and
only if $X_{ab}=X_a+\Ad(g_a)da(X_b)$ for all $a,b\in A$. Since
$\Ad(g_a)da$ is the differential of the action of $a$ on
$G_\alpha$, the $A$-module obtained by twisting $G$ using
$\alpha$, $a\mapsto\Ad(g_a)da$ is a linear action of $A$ on $\g$.
Hence $\ker((d\varphi)_\alpha)=Z^1(A,\g)$.

For $\alpha\in Z^1(A,G)$ in the previous paragraph, we define
$\psi:G\rightarrow\Map(A,G)$ by $\psi(g)(a)=g^{-1}g_aa(g)$. Then
$\psi(e)=\alpha$ and $\psi(g)\in Z^1(A,G)$ for all $g\in G$. For
$X\in\g$, we have, under the identification
$T_\alpha\Map(A,G)\cong\Map(A,\g)$,
\begin{align*}
&(d\psi)_e(X)(a)\\
=&\frac{d}{dt}\Big|_{t=0}e^{-tX}g_aa(e^{tX})g_a^{-1}\\
=&\frac{d}{dt}\Big|_{t=0}e^{-tX}e^{t\Ad(g_a)da(X)}\\
=&-X+\Ad(g_a)da(X).
\end{align*}
Hence $\mathrm{Im}((d\psi)_e)=B^1(A,\g)$, the set of
$1$-coboundaries with respect to the linear action of $A$ on $\g$
mentioned above.

Note that $H^1(A,\g)=Z^1(A,\g)/B^1(A,\g)$ is a linear space over
$\RR$. But by \cite{Se1}, Chapter VIII, Section 2, Corollary 1,
all elements of $H^1(A,\g)$ are annihilated by $|A|$, the order of
$A$. So $H^1(A,\g)=0$, that is $Z^1(A,\g)=B^1(A,\g)$, and then
$\ker((d\varphi)_\alpha)=\mathrm{Im}((d\psi)_e)$. By Weil's Lemma,
every connected component of $Z^1(A,G)=\varphi^{-1}(e)$ is a
closed submanifold of $\Map(A,G)$.

Now we prove that every cohomology class in $Z^1(A,G)$ is a
connected component of $Z^1(A,G)$. In fact, for $\alpha\in
Z^1(A,G)$, the cohomology class of $\alpha$ is the image of the
map $\psi$ defined above. By Weil's Lemma, there exist a
neighborhood $U$ of $\alpha$ in $\Map(A,G)$ and a neighborhood $V$
of $e$ in $G$ such that $\psi(V)=U\cap Z^1(A,G)$. So the
cohomology class of $\alpha$ is open in $Z^1(A,G)$. Since the
cohomology class of $\alpha$ is the complement of the union of
other cohomology classes, which are all open in $Z^1(A,G)$, it is
also closed in $Z^1(A,G)$. But as the image of $\psi$, it is
connected. Hence it is a connected component of $Z^1(A,G)$.

Finally we prove that $Z^1(A,G)$ has only finitely many connected
components. Since we have proved that the cardinality of the set
of connected components of $Z^1(A,G)$ coincides with that of
$H^1(A,G)$, by Theorem \ref{T:A-bi}, it is sufficient to prove the
case for which $G$ is compact. But if $G$ is compact, so is
$\Map(A,G)$, and as we have showed above, every connected
component of the closed subset $Z^1(A,G)$ of $\Map(A,G)$ is open
in $Z^1(A,G)$. This forces that $Z^1(A,G)$ has only finitely many
connected components. The proof of the theorem is finished.
\end{proof}

\begin{corollary}\label{C:A-finite}
Let $A$ be a finite group, and let $G$ be a connected Lie group
with an $A$-module structure. Then $H^1(A,G)$ is finite.\qed
\end{corollary}

\begin{remark}
Professor Jiu-Kang Yu had a simple algebro-geometric proof of
Theorem \ref{A-cocycle} when $G$ is compact.
\end{remark}

%%%%%%%%%%%%%%%%%%%%%%%%%%%%%%%%%%%%%%%%%%%%%%%%%%%%%%%  Section 4

\section{Nonabelian cohomology of $\ZZ$ with coefficients in Lie
groups}

Let $V$ be a finite-dimensional real vector space, $L$ a linear
endomorphism of $V$. For a complex number $\lambda$, $L$ is
\emph{$\lambda$-semisimple} if $\ker(\lambda I-L)=\ker((\lambda
I-L)^2)$ in the complexification $V_\mathbb{C}$ of $V$. Note that
$L$ is semisimple if and only if it is $\lambda$-semisimple for
every $\lambda\in\mathbb{C}$. For a connected Lie group $G$ and a
complex number $\lambda$, an automorphism $\sigma$ of $G$ is
\emph{$\lambda$-semisimple} (\emph{semisimple}, resp.) if the
linear endomorphism $d\sigma$ on the Lie algebra $\g$ of $G$ is
$\lambda$-semisimple (semisimple, resp.). Suppose $A$ is a cyclic
group and $G$ is a connected Lie group with an $A$-module
structure. The action of $A$ on $G$ is \emph{$1$-semisimple}
(\emph{semisimple}, resp.) if the action of a generator of $A$ on
$G$ is $1$-semisimple (semisimple, resp.). Note that this
definition is independent of the choice of the generator of $A$.

All automorphisms of a connected compact semisimple Lie group are
$1$-semisimple, due to the fact that the automorphism group of a
compact semisimple Lie group is compact. On the other hand, not
all automorphisms of a compact torus of dimension $\geq2$ are
$1$-semisimple. For example, the automorphism of $\RR^2/\ZZ^2$
induced by ${\SMALL\begin{pmatrix}
1 & 1\\
0 & 1
\end{pmatrix}}$ is not $1$-semisimple. In fact, For a connected Lie
group $G$ with an automorphism $\sigma$, $\sigma$ is
$1$-semisimple if and only if the restriction of $\sigma$ to the
connected component of the center of $G$ is $1$-semisimple.

The main content of this section is to prove the following
theorem.

\begin{theorem}\label{T:Z-sur}
Let $G$ be a connected compact Lie group with a $\ZZ$-module
structure, $T$ a maximal torus of $G_0^\ZZ$. Let
$i_1:H^1(\ZZ,T)\rightarrow H^1(\ZZ,G)$ be the canonical map. Then
the
following statements are equivalent.\\
(i) The $\ZZ$-action on $G$ is $1$-semisimple;\\
(ii) The map $i_1$ is surjective;\\
(iii) The kernel of $i_1$ is finite;\\
(iv) All fibers of $i_1$ are finite.
\end{theorem}

To prove Theorem \ref{T:Z-sur}, we need several lemmas. The first
lemma is a key step, whose proof needs the following facts.

\begin{fact}[\cite{Ko}, Chapter II, Theorem 5.1]\label{F:fix-geodesic}
Let $M$ be a Riemannian manifold with an isometry $f$. Then every
connected component of $\mathrm{Fix}(f)$ is a closed totally
geodesic submanifold of $M$, where $\mathrm{Fix}(f)$ is the set of
fixed points of $f$.
\end{fact}

\begin{fact}[\cite{Ko}, Chapter II, Section 5, page 63]\label{F:Lefschetz}
Let $M$ be a compact Riemannian manifold with an isometry $f$.
Then $\mathrm{Lef}(f)=\chi(\mathrm{Fix}(f))$, where
$\mathrm{Lef}(f)$ is the Lefschetz number of $f$,
$\chi(\mathrm{Fix}(f))$ is the Euler number of $\mathrm{Fix}(f)$.
\end{fact}

\begin{fact}[\cite{KN}, Chapter VII, Corollary 8.10]\label{F:geodesic-homogeneous}
Let $M$ be a compact homogeneous Riemannian manifold. Then every
closed totally geodesic submanifold of $M$ is a compact
homogeneous Riemannian manifold.
\end{fact}

\begin{fact}[\cite{Sa}, Theorem II]\label{F:Euler}
Let $N$ be a connected compact homogeneous Riemannian manifold.
Then $\chi(N)\geq0$. Write $N$ as the form $N=K/H$, where $K$ and
$H$ are compact Lie groups, $K$ is connected. Then $\chi(N)>0$ if
and only if $\mathrm{rank}(K)=\mathrm{rank}(H_0)$.
\end{fact}

\begin{lemma}\label{L:ss-sur}
Let $G$ be a connected compact semisimple Lie group with a
$\ZZ$-module structure, $T$ a maximal torus of $G_0^\ZZ$. Then the
canonical map $H^1(\ZZ,T)\rightarrow H^1(\ZZ,G)$ is surjective.
\end{lemma}

\begin{sublemma}\label{L:isomorphism}
Suppose $G$ is a connected compact semisimple Lie group with a
closed subgroup $H$. Denote $\Aut(G,H)=\{\theta\in
\Aut(G)|\theta(H)=H\}$, which acts on $G/H$ by
$\rho_\theta(gH)=\theta(g)H$. Then there is an $\Ad(G)$-invariant
inner product on $\g$ such that the induced $G$-invariant
Riemannian structure on $G/H$ is also $\Aut(G,H)$-invariant.
\end{sublemma}

\begin{proof}
We construct the semidirect product $G\rtimes \Aut(G,H)$, and
define the map $\varphi:(G\rtimes \Aut(G,H))\times G/H\rightarrow
G/H$ by $\varphi((g,\theta),g'H)=g\theta(g')H$. We claim that
$\varphi$ is an action of $G\rtimes \Aut(G,H)$ on $G/H$. In fact,
if we denote $\varphi_{(g,\theta)}=\varphi((g,\theta),\cdot)$,
then
\begin{align*}
&\varphi_{(g_1,\theta_1)}\circ\varphi_{(g_2,\theta_2)}(g'H)\\
=&\varphi_{(g_1,\theta_1)}(g_2\theta_2(g')H)\\
=&g_1\theta_1(g_2)\theta_1(\theta_2(g'))H\\
=&\varphi_{(g_1\theta_1(g_2),\theta_1\circ\theta_2)}(g'H)\\
=&\varphi_{(g_1,\theta_1)(g_2,\theta_2)}(g'H).
\end{align*}
This verifies the claim. Note that the action of $G\rtimes
\Aut(G,H)$ on $G/H$ is transitive, and the isotropic subgroup
associated with the point $H\in G/H$ is $H\rtimes \Aut(G,H)$. So
there is a natural isomorphism $G/H\cong G\rtimes
\Aut(G,H)/H\rtimes \Aut(G,H)$. Since $G$ is compact semisimple,
$\Aut(G)$ is compact. As a closed subgroup of $\Aut(G)$,
$\Aut(G,H)$ is also compact. So $G\rtimes \Aut(G,H)$ is compact.
Choose an $\Ad(G\rtimes \Aut(G,H))$-invariant inner product
$B(\cdot,\cdot)$ on the Lie algebra $\mathrm{L}(G\rtimes
\Aut(G,H))$ of $G\rtimes \Aut(G,H)$. Then the restriction of $B$
on $\mathrm{L}(H\rtimes \Aut(G,H))^\perp$ induces a $G\rtimes
\Aut(G,H)$-invariant Riemannian structure on $G/H$. Note that if
we identify $G\times\{id\}$ with $G$ and $\{e\}\times \Aut(G,H)$
with $\Aut(G,H)$, the restrictions of the $G\rtimes
\Aut(G,H)$-action to $G$ and $\Aut(G,H)$ coincide with the natural
actions of $G$ and $\Aut(G,H)$ on $G/H$. So the $G\rtimes
\Aut(G,H)$-invariant Riemannian structure on $G/H$ constructed
above is $G$-invariant and $\Aut(G,H)$-invariant.

Note that the restriction of the inner product $B$ on the Lie
algebra $\g$ of $G$ is $\Ad(G)$-invariant, and the restriction of
$B$ on $\ah^\perp$ induces a $G$-invariant Riemannian structure on
$G/H$ which is same as that was constructed above. This proves the
sublemma.
\end{proof}

\begin{sublemma}\label{L:F_0}
Under the same conditions as in Sublemma \ref{L:isomorphism}, if
$\theta\in \Aut(G,H)$ satisfies $H\subset G^\theta$, then the
connected component of $\mathrm{Fix}(\rho_\theta)\subset G/H$
containing the point $H$ is $G^\theta_0/H$.
\end{sublemma}

\begin{proof}
Endow $G/H$ with a $G$-invariant and $\Aut(G,H)$-invariant
Riemannian structure which is induced by an $\Ad(G)$-invariant
inner product on $\g$. By Corollary 3.6, Theorem 2.10, and
Corollary 2.5 in Chapter X of \cite{KN}, every geodesic of $G/H$
starting from $H$ is of the form $e^{tX}H$ for some $X\in\g$.

Denote the connected component of $\mathrm{Fix}(\rho_\theta)$
containing $H$ by $F_0$. By Fact \ref{F:fix-geodesic}, $F_0$ is a
closed totally geodesic submanifold of $G/H$. Suppose $gH\in F_0$.
Let $e^{tX}H$ be a geodesic through $H$ and $gH$ such that
$e^{tX}H\in F_0, \forall t\in\mathbb{R}$. Then
$e^{tX}H=\rho_\theta(e^{tX}H)=\theta(e^{tX})H=e^{td\theta(X)}H$.
This implies $e^{-tX}e^{td\theta(X)}\in H$. So
$(d\theta-1)(X)\in\ah\subset\mathrm{L}(G^\theta)=\ker(d\theta-1)$.
Since $G$ is compact semisimple, $d\theta$ is semisimple. This
implies $X\in\mathrm{L}(G^\theta)$, and then $gH\in G^\theta_0/H$.
So $F_0\subset G^\theta_0/H$. The inverse direction is obvious.
\end{proof}

\begin{proof}[Proof of Lemma \ref{L:ss-sur}]
Let $\sigma$ be a generator of $\ZZ$. Then under the
identification $Z^1(\ZZ,G)\cong G$ (see Section 2), two cocycles
$z,z'\in G$ are cohomologous if there exists $g\in G$ such that
$z'=g^{-1}z\sigma(g)$. So to prove the lemma, it is sufficient to
show that for each $z\in G$, there exists $g\in G$ such that
$g^{-1}z\sigma(g)\in T$.

Let $f$ be the diffeomorphism of $G/T$ defined by
$f(gT)=\sigma(g)T$. Since $\sigma\in \Aut(G,T)$ and $T\subset
G^\ZZ$, by Sublemma \ref{L:isomorphism} and Sublemma \ref{L:F_0},
there is a $G$-invariant Riemannian structure on $G/T$ which is
also $f$-invariant, and the connected component of
$\mathrm{Fix}(f)$ containing $T$ is $G^\ZZ_0/T$. Let
$\{F_0,F_1,\cdots\}$ be the set of all the connected components of
$\mathrm{Fix}(f)$ with $F_0=G^\ZZ_0/T$. Since
$\mathrm{rank}(G^\ZZ_0)=\mathrm{rank}(T)$, by Fact \ref{F:Euler},
we have $\chi(G^\ZZ_0/T)>0$. Since $f$ is an isometry with respect
to the above-mentioned Riemannian structure on $G/T$, by Fact
\ref{F:fix-geodesic}, each $F_i$ is a closed totally geodesic
submanifold of $G/T$. Then by Fact \ref{F:geodesic-homogeneous},
each $F_i$ is a compact homogeneous Riemannian manifold. So
$\chi(F_i)\geq0$ for each $i$, by Fact \ref{F:Euler}. Applying the
formula in Fact \ref{F:Lefschetz}, we have
$$\mathrm{Lef}(f)=\chi(\mathrm{Fix}(f))=\sum_{i\geq0}\chi(F_i)
=\chi(G^\ZZ_0/T)+\sum_{i>0}\chi(F_i)\geq\chi(G^\ZZ_0/T)>0.$$

For each $z\in G$, define the diffeomorphism $f_z$ of $G/T$ by
$f_z(gT)=z\sigma(g)T$. Since $G$ is connected, $f_z$ is homotopic
to $f=f_e$, and then $\mathrm{Lef}(f_z)=\mathrm{Lef}(f)>0$. By the
Lefschetz Fixed Point Theorem, $f_z$ has a fixed point $gT$, that
is, $z\sigma(g)T=gT$. This means that $g^{-1}z\sigma(g)\in T$. The
lemma is proved.
\end{proof}

\begin{lemma}\label{L:torus-sur}
Let $G$ be a compact torus with a $\ZZ$-module structure,
$T=G_0^\ZZ$.\\
(i) If the $\ZZ$-action on $G$ is $1$-semisimple, then
$H^1(\ZZ,T)\rightarrow H^1(\ZZ,G)$ is surjective.\\
(ii) If there exists a discrete $\ZZ$-submodule $\Gamma$ of $G$
such that $H^1(\ZZ,T\Gamma)\rightarrow H^1(\ZZ,G)$ is surjective,
then the $\ZZ$-action on $G$ is $1$-semisimple.
\end{lemma}

\begin{proof}
(i) We also choose a generator $\sigma$ of $\ZZ$ and identify
$Z^1(\ZZ,G)$ with $G$. Let $z\in G$ be a cocycle, and choose an
$X\in\g$ such that $e^X=z$. Since $\sigma$ is $1$-semisimple,
$\ker(1-d\sigma)\oplus\mathrm{Im}(1-d\sigma)=\g$. So we can write
$X$ as $X=Y+(1-d\sigma)(Z)$, that is, $-Z+X+d\sigma(Z)=Y$, where
$Y\in\ker(1-d\sigma), Z\in\g$. Since $G$ is abelian, this implies
that $e^{-Z}z\sigma(e^Z)=e^Y\in T$. This proves (i).

(ii) We prove it by contradiction. Suppose that the $\ZZ$-action
on $G$ is not $1$-semisimple. So
$\ker(1-d\sigma)+\mathrm{Im}(1-d\sigma)\subsetneqq\g$. Then for
each discrete $\ZZ$-submodule $\Gamma$ of $G$,
$\ker(1-d\sigma)+\mathrm{Im}(1-d\sigma)+\exp^{-1}(\Gamma)\subsetneqq\g$.
Choose an $X\in\g$ which does not belong to the left hand side of
the above equation. Then it is easy to check that for every $g\in
G$, $g^{-1}e^X\sigma(g)\notin T\Gamma$, that is,
$H^1(\ZZ,T\Gamma)\rightarrow H^1(\ZZ,G)$ is not surjective. This
conflicts with the conditions, and then (ii) is proved.
\end{proof}

\begin{lemma}\label{L:cover-sur}
Let $G,G'$ be connected Lie groups with $\ZZ$-module structures,
and let $\pi:G'\rightarrow G$ be a
$\ZZ$-epimorphism.\\
(i) If $H'$ is a submodule of $G'$ such that
$H^1(\ZZ,H')\rightarrow H^1(\ZZ,G')$ is surjective, then
$H^1(\ZZ,\pi(H'))\rightarrow H^1(\ZZ,G)$ is surjective.\\
(ii) If $H$ is a submodule of $G$ such that $H^1(\ZZ,H)\rightarrow
H^1(\ZZ,G)$ is surjective, then $H^1(\ZZ,\pi^{-1}(H))\rightarrow
H^1(\ZZ,G')$ is surjective.
\end{lemma}

\begin{proof}
The first assertion follows from the following commutative
diagram, in which the two columns are surjective.
$$
\xymatrix{ H^1(\ZZ,H')\ar[d]\ar[r]
& H^1(\ZZ,G')\ar[d]\\
H^1(\ZZ,\pi(H')) \ar[r] & H^1(\ZZ,G)}
$$

To prove the second assertion, let $z'\in G'$ be a cocycle. If
$H^1(\ZZ,H)\rightarrow H^1(\ZZ,G)$ is surjective, then there
exists $g\in G$ such that $g^{-1}\pi(z')\sigma(g)\in H$. Choose a
$g'\in G'$ with $\pi(g')=g$, then
$\pi(g'^{-1}z'\sigma(g'))=g^{-1}\pi(z')\sigma(g)\in H$, that is,
$g'^{-1}z'\sigma(g')\in\pi^{-1}(H)$. Hence
$H^1(\ZZ,\pi^{-1}(H))\rightarrow H^1(\ZZ,G')$ is surjective.
\end{proof}

\begin{lemma}\label{L:ss-finite}
Let $G$ be a connected compact semisimple Lie group with a
$\ZZ$-module structure, $T$ a maximal torus of $G_0^\ZZ$. Then all
fibers of $H^1(\ZZ,T)\rightarrow H^1(\ZZ,G)$ are finite.
\end{lemma}

\begin{sublemma}\label{L:manifold}
Let $M$ be a compact smooth manifold with two closed smooth
submanifolds $N_1,N_2$. If for each $p\in N_1\cap N_2$,
$T_pN_1\cap T_pN_2=0$, then $N_1\cap N_2$ is a finite set.
\end{sublemma}

\begin{proof}
Let $p\in N_1\cap N_2$. Choose a coordinate chart
$(U,x_1,\cdots,x_n)$ of $M$ containing $p$ such that $x_i(p)=0,
1\leq i\leq n$, and such that $U\cap N_1=\{q\in
U|x_{n_1+1}(q)=\cdots=x_n(q)=0\}$, where $n=\dim M, n_1=\dim N_1$.
Define $f:U\rightarrow \RR^{n-n_1}$ by
$f(q)=(x_{n_1+1}(q),\cdots,x_n(q))$. Then $T_pN_1\cap T_pN_2=0$
implies that $f|_{U\cap N_2}:U\cap N_2\rightarrow \RR^{n-n_1}$ is
an immersion at $p$. Hence there is an open neighborhood $V\subset
U\cap N_2$ of $p$ in $N_2$ such that $f|_V$ is injective. In
particular, $V\cap N_1=\{p\}$. This proves that $N_1\cap N_2$ is
discrete in $N_2$. Since $N_2$ is compact, $N_1\cap N_2$ is
finite.
\end{proof}

\begin{proof}[Proof of Lemma \ref{L:ss-finite}]
A fiber of $H^1(\ZZ,T)\rightarrow H^1(\ZZ,G)$ is the intersection
of $T$ with a cohomology class in $Z^1(\ZZ,G)$. But a cohomology
class in $Z^1(\ZZ,G)$ is an orbit of the twisted conjugate action
$\tau$ associated with $\sigma$ (see Section 2), which is a closed
submanifold. So by Sublemma \ref{L:manifold}, it is sufficient to
prove that for each $t\in T$, $T_tO_t\cap T_tT=0$, where $O_t$ is
the orbit of the twisted conjugate action containing $t$.

Let $l_t$ be the left translation on $G$ induced by $t$. Since
\begin{align*}
&\frac{d}{ds}\Big|_{s=0}\tau_{e^{sX}}(t)\\
=&(dl_t)_e\left(\frac{d}{ds}\Big|_{s=0}t^{-1}e^{sX}t\sigma(e^{sX})^{-1}\right)\\
=&(dl_t)_e\left(\frac{d}{ds}\Big|_{s=0}e^{s\Ad(t^{-1})X}e^{-sd\sigma(X)}\right)\\
=&(dl_t)_e(\Ad(t^{-1})-d\sigma)(X),
\end{align*}
we have $T_tO_t\cap
T_tT=(dl_t)_e(\textrm{Im}(\Ad(t^{-1})-d\sigma)\cap\at)$, and it is
sufficient to show that
$\textrm{Im}(\Ad(t^{-1})-d\sigma)\cap\at=0$. Since $G$ is compact
semisimple, $\Ad(t^{-1})$ and $d\sigma$ are semisimple. The fact
$\sigma(t)=t$ implies that $\Ad(t^{-1})$ and $d\sigma$ commute, so
$\Ad(t^{-1})-d\sigma$ is also semisimple. But
$\at\subset\ker(\Ad(t^{-1})-d\sigma)$, so
$\textrm{Im}(\Ad(t^{-1})-d\sigma)\cap\at=0$. This proves the
lemma.
\end{proof}

\begin{lemma}\label{L:torus-finite}
Let $G$ be a compact torus with a $\ZZ$-module structure,
$T=G_0^\ZZ$. Then the kernel of $i_1:H^1(\ZZ,T)\rightarrow
H^1(\ZZ,G)$ is finite if and only if the $\ZZ$-action on $G$ is
$1$-semisimple.
\end{lemma}

\begin{proof}
Since $G$ is abelian, we have the following exact sequence of
cohomology \emph{groups}:
$$
0\rightarrow H^0(\ZZ,T)\stackrel{i_0}\rightarrow
H^0(\ZZ,G)\rightarrow H^0(\ZZ,G/T)\stackrel{\delta}\rightarrow
H^1(\ZZ,T)\stackrel{i_1}\rightarrow H^1(\ZZ,G).
$$
Since $\ker(\delta)\cong\mathrm{coker}(i_0)=G^\ZZ/T$ is finite,
\begin{align*}
&\text{$\ker(i_1)
=\mathrm{Im}(\delta)$ is finite}\\
\Leftrightarrow & \text{$H^0(\ZZ,G/T)$ is finite}\\
\Leftrightarrow & \text{the induced linear automorphism of
$d\sigma$ on $\g/\at$ has no eigenvalue $1$}\\
\Leftrightarrow & \text{$d\sigma$ is $1$-semisimple}\\
\Leftrightarrow & \text{the $\ZZ$-action on $G$ is
$1$-semisimple.}
\end{align*}
\end{proof}

\begin{proof}[Proof of Theorem \ref{T:Z-sur}]
We follow the line
``(i)$\Rightarrow$(ii)$\Rightarrow$(i)$\Rightarrow$(iv)
$\Rightarrow$(iii)$\Rightarrow$(i)".

Let $\g_s=[\g,\g]$, $\g_t$ be the center of $\g$. Then
$\g=\g_s\oplus\g_t$. Let $G_s, G_t$ be the connected subgroups of
$G$ with Lie algebras $\g_s$ and $\g_t$, respectively. Then $G_s$
is semisimple, $G_t$ is a compact torus, and it is obvious that
$G_s$ and $G_t$ are $\ZZ$-submodules of $G$. Let $G'=G_s\times
G_t$ be the direct product of $\ZZ$-modules, and define
$\pi:G'\rightarrow G$ by $\pi(g_s,g_t)=g_sg_t$. Then $\pi$ is a
$\ZZ$-epimorphism which is also a finite covering homomorphism of
Lie groups. Let $T_s=(T\cap G_s)_0, T_t=(T\cap G_t)_0,
T'=T_s\times T_t$, then $T_s, T_t$, and $T'$ are maximal tori of
$(G_s^\ZZ)_0,(G_t^\ZZ)_0$, and $G'^\ZZ_0$, respectively, and we
have $\pi(T')=T$. Note that the $\ZZ$-action on $G$ is
$1$-semisimple if and only if the $\ZZ$-action on $G_t$ is
$1$-semisimple.

``(i)$\Rightarrow$(ii)". Suppose that the $\ZZ$-action on $G$ is
$1$-semisimple. By Lemma \ref{L:ss-sur} and part (i) of Lemma
\ref{L:torus-sur}, $H^1(\ZZ,T_s)\rightarrow H^1(\ZZ,G_s)$ and
$H^1(\ZZ,T_t)\rightarrow H^1(\ZZ,G_t)$ are surjective. From the
commutative diagram
\begin{equation}\label{E:Gs*Gt}
\xymatrix{ H^1(\ZZ,T')\ar[d]\ar[r]^\cong
& H^1(\ZZ,T_s)\ar[d] & \times & H^1(\ZZ,T_t)\ar[d]\\
H^1(\ZZ,G')\ar[r]^\cong & H^1(\ZZ,G_s) & \times & H^1(\ZZ,G_t),}
\end{equation}
we know that $H^1(\ZZ,T')\rightarrow H^1(\ZZ,G')$ is surjective.
By part (i) of Lemma \ref{L:cover-sur}, $H^1(\ZZ,T)\rightarrow
H^1(\ZZ,G)$ is surjective.

``(ii)$\Rightarrow$(i)". Suppose that $H^1(\ZZ,T)\rightarrow
H^1(\ZZ,G)$ is surjective. By part (ii) of Lemma
\ref{L:cover-sur}, $H^1(\ZZ,T'\cdot\ker(\pi))\rightarrow
H^1(\ZZ,G')$ is surjective. Applying part (i) of Lemma
\ref{L:cover-sur} to the projection $p_2:G'=G_s\times
G_t\rightarrow G_t$, we know that $H^1(\ZZ,T_t\cdot
p_2(\ker(\pi)))\rightarrow H^1(\ZZ,G_t)$ is surjective. But
$p_2(\ker(\pi))$ is a discrete $\ZZ$-submodule of $G_t$, by part
(ii) of Lemma \ref{L:torus-sur}, the $\ZZ$-action on $G_t$ is
$1$-semisimple. Hence the $\ZZ$-action on $G$ is $1$-semisimple.

``(i)$\Rightarrow$(iv)". Suppose that the $\ZZ$-action on $G$ is
$1$-semisimple. We show that all fibers of $H^1(\ZZ,T)\rightarrow
H^1(\ZZ,G)$ are finite. First we prove the kernel of this map is
finite. By Lemma \ref{L:ss-finite}, Lemma \ref{L:torus-finite},
and the commutative diagram \eqref{E:Gs*Gt}, we know that the
kernel of $H^1(\ZZ,T')\rightarrow H^1(\ZZ,G')$ is finite. Since
$\pi:G'\rightarrow G$ is a finite cover, $\ker(\pi)$ and
$\ker(\pi|_{T'})$ are finite. Hence $H^1(\ZZ,\ker(\pi))$ and
$H^1(\ZZ,\ker(\pi|_{T'}))$ are finite. By the commutative diagram
\begin{equation}\label{E:finite}
\xymatrix{ H^1(\ZZ,\ker(\pi|_{T'}))\ar[d] &
H^1(\ZZ,\ker(\pi))\ar[d]\\
H^1(\ZZ,T')\ar[d]\ar[r] & H^1(\ZZ,G')\ar[d]\\
H^1(\ZZ,T)\ar[r] & H^1(\ZZ,G),}
\end{equation}
in which the two columns are exact, we know that the kernel of
$H^1(\ZZ,T)\rightarrow H^1(\ZZ,G)$ is finite.

Now let $[\alpha]\in H^1(\ZZ,G)$. By the step
``(i)$\Rightarrow$(ii)" above, we know the fiber
$i_1^{-1}([\alpha])$ is not empty. Choose a $[\beta]\in
i_1^{-1}([\alpha])$, where $\beta\in[\beta]$. Then the twisting
using $\beta$ transform the fiber $i_1^{-1}([\alpha])$ to the
kernel of $H^1(\ZZ,T_\beta)\rightarrow H^1(\ZZ,G_{i\circ\beta})$
(see Section 2). Since $T_\beta\subset(G_{i\circ\beta}^\ZZ)_0$, we
can choose a maximal torus $T'_\beta$ of $(G_{i\circ\beta}^\ZZ)_0$
with $T_\beta\subset T'_\beta$. It is easy to see that this
twisting process does not change the $1$-semisimplicity of the
$\ZZ$-action on $G$, so the kernel of
$H^1(\ZZ,T'_\beta)\rightarrow H^1(\ZZ,G_{i\circ\beta})$ is finite.
Note that $H^1(\ZZ,T_\beta)\rightarrow H^1(\ZZ,T'_\beta)$ is
injective, hence the kernel of $H^1(\ZZ,T_\beta)\rightarrow
H^1(\ZZ,G_{i\circ\beta})$ is finite. This proves that
$i_1^{-1}([\alpha])$ is finite.

``(iv)$\Rightarrow$(iii)". Obvious.

``(iii)$\Rightarrow$(i)". Suppose that the kernel of
$H^1(\ZZ,T)\rightarrow H^1(\ZZ,G)$ is finite. By the commutative
diagram \eqref{E:finite}, we know the kernel of
$H^1(\ZZ,T')\rightarrow H^1(\ZZ,G')$ is finite. And by the
commutative diagram \eqref{E:Gs*Gt}, we know the kernel of
$H^1(\ZZ,T_t)\rightarrow H^1(\ZZ,G_t)$ is also finite. By Lemma
\ref{L:torus-finite}, the $\ZZ$-action on $G_t$ is $1$-semisimple.
Hence the $\ZZ$-action on $G$ is $1$-semisimple.
\end{proof}

Theorem \ref{T:Z-sur} has the following corollaries.

\begin{corollary}\label{C:semisimple}
Let $G$ be a connected compact semisimple Lie group with a
$\ZZ$-module structure, $T$ a maximal torus of $G_0^\ZZ$. Then the
canonical map $i_1:H^1(\ZZ,T)\rightarrow H^1(\ZZ,G)$ is surjective
with finite fibers.
\end{corollary}

\begin{proof}
Since $G$ is compact semisimple, the automorphism group $\Aut(G)$
of $G$ is compact. So the action of each $\sigma\in\Aut(G)$ on $G$
is semisimple. As a sub-action of $\Aut(G)$, the $\ZZ$-action on
$G$ is semisimple, hence $1$-semisimple, and the corollary follows
from Theorem \ref{T:Z-sur} immediately.
\end{proof}

\begin{corollary}\label{C:Z-coboundary}
Let $G$ be a connected compact Lie group with a $\ZZ$-module
structure, $T$ a maximal torus of $G_0^\ZZ$. Suppose the
$\ZZ$-action on $G$ is $1$-semisimple. Then every connected
component of $H^0(\ZZ,G/T)$ is a closed homogeneous submanifold of
$G/T$, and the restriction of the coboundary operator
$\delta:H^0(\ZZ,G/T)\rightarrow H^1(\ZZ,T)$ to each connected
component of $H^0(\ZZ,G/T)$ is constant.
\end{corollary}

\begin{proof}
From the exact sequence
$$
H^0(\ZZ,G/T)\stackrel{\delta}\rightarrow
H^1(\ZZ,T)\stackrel{i_1}\rightarrow H^1(\ZZ,G)
$$
we know that the image of $\delta$ coincides with the kernel of
$i_1$, which is finite by Theorem \ref{T:Z-sur}. So the image of
each connected component of $H^0(\ZZ,G/T)$ is finite. But under
the identification $H^1(\ZZ,T)$ with $T$, $\delta$ is continuous.
This forces that the restriction of $\delta$ to each connected
component of $H^0(\ZZ,G/T)$ is constant.

Let $M$ be a connected component of $H^0(\ZZ,G/T)$, and suppose
that $\delta(M)=\{t\}\subset T\cong H^1(\ZZ,T)$. We compute what
the subset $\delta^{-1}(t)$ of $G/T$ is. Choose a point $g_0T\in
M$. Since $\delta:H^0(\ZZ,G/T)\rightarrow H^1(\ZZ,T)\cong T$ has
the expression $\delta(gT)=g^{-1}\sigma(g)$, we have
$g_0^{-1}\sigma(g_0)=t$. Now suppose $gT\in H^0(\ZZ,G/T)$, then
\begin{align*}
&\delta(gT)=g^{-1}\sigma(g)=t\\
\Leftrightarrow &
(g_0^{-1}g)^{-1}g_0^{-1}\sigma(g_0)\sigma(g_0^{-1}g)=t\\
\Leftrightarrow & t\sigma(g_0^{-1}g)t^{-1}=g_0^{-1}g.
\end{align*}
Since $H=\{h\in G|t\sigma(h)t^{-1}=h\}$ is a closed subgroup of
$G$ containing $T$, $\delta^{-1}(t)$ is diffeomorphic to $H/T$.
Hence $M$, as a connected component of $\delta^{-1}(t)$, is
diffeomorphic to $H_0/T$, which is a closed homogeneous
submanifold of $G/T$. This proves the corollary.
\end{proof}

%%%%%%%%%%%%%%%%%%%%%%%%%%%%%%%%%%%%%%%%%%%%%%%%%%%%%%%  Section 5

\section{Nonabelian cohomology of $\Zn$ with coefficients in Lie
groups}

Let $A=\Zn$, $G$ be a connected Lie group with an $A$-module
structure. Then all assertions which have been proved in Section 3
hold for $H^1(A,G)$. In this section we prove further properties
of $H^1(\Zn,G)$. We first introduce a notion.

Let $H$ be a connected Lie group. A subgroup $T$ of $H$ is a
\emph{maximal compact torus} of $H$ if $T$ is a compact torus and
there is no other compact torus $T'$ of $H$ such that
$T\subsetneqq T'$. Note that if $H$ is compact, this notion
coincides with the commonly used notion of maximal torus.

It is obvious that maximal compact tori always exist. We claim
that they are unique up to conjugation.

\begin{proposition}\label{P:torus}
Let $H$ be a connected Lie group, $T$ a maximal compact torus of
$H$. Then for any compact torus $T'$ of $H$, there exists an $h\in
H$ such that $hT'h^{-1}\subset T$. In particular, any two maximal
compact tori of $H$ are conjugate.
\end{proposition}

\begin{proof}
Let $K$ be a maximal compact subgroup of $G$ containing $T$. By
Fact \ref{F:max-cpt}, there is an $h_1\in H$ such that
$h_1T'h_1^{-1}\subset K$. Since $K$ is connected compact and $T$
is a maximal torus of $K$, there is a $k\in K$ such that
$(kh_1)T'(kh_1)^{-1}\subset T$.
\end{proof}

\begin{theorem}\label{T:Zn-sur}
Let $G$ be a connected Lie group with a $\Zn$-module structure,
$T$ a maximal compact torus of $G_0^{\Zn}$. Then the canonical map
$i_1:H^1(\Zn,T)\rightarrow H^1(\Zn,G)$ is surjective.
\end{theorem}

\begin{proof}
Let $\sigma$ be a generator of $\Zn$. Then under the
identification $Z^1(\Zn,G)\cong Z=\{g\in
G|g\sigma(g)\cdots\sigma^{n-1}(g)=e\}$ (see Section 2), two
cocycles $z,z'\in Z$ are cohomologous if there exists $g\in G$
such that $z'=g^{-1}z\sigma(g)$. So to prove the theorem, it is
sufficient to show that for each $z\in Z$, there exists $g\in G$
such that $g^{-1}z\sigma(g)\in T$.

Let $z\in Z$. By Theorem \ref{T:A-bi}, there exists a maximal
compact subgroup $K$ of $G$ which is also a $\Zn$-submodule of
$G$, and there also exists $g_1\in G$ such that
$g_1^{-1}z\sigma(g_1)\in K$.

Denote $z_1=g_1^{-1}z\sigma(g_1)$. Let $T'$ be a maximal compact
torus of $K_0^{\Zn}$. Since $\sigma$ is of finite order, it is
semisimple, hence $1$-semisimple. By Theorem \ref{T:Z-sur}, there
is a $k\in K$ such that $k^{-1}z_1\sigma(k)\in T'$.

Denote $z_2=k^{-1}z_1\sigma(k)$. Since $T'$ is a compact torus of
$G_0^{\Zn}$, by Proposition \ref{P:torus}, there exists a $g_2\in
G_0^{\Zn}$ such that $g_2^{-1}T'g_2\subset T$. In particular,
$g_2^{-1}z_2g_2\in T$. Denote $g=g_1kg_2$, and notice that
$\sigma(g_2)=g_2$, we have $g^{-1}z\sigma(g)\in T$. This proves
the theorem.
\end{proof}

\begin{corollary}\label{C:Zn-trivial}
Let $G$ be a connected and simply connected solvable Lie group
with a $\Zn$-module structure. Then $H^1(\Zn,G)$ is trivial.
\end{corollary}

\begin{proof}
Let $T$ be a maximal compact torus of $G_0^{\Zn}$. By Theorem
\ref{T:Zn-sur}, it is sufficient to show that $H^1(\Zn,T)$ is
trivial. But a connected and simply connected solvable Lie group
has no nontrivial compact subgroup (\cite{Ho}, Chapter XII,
Theorem 2.3). So $T$ is trivial, and hence $H^1(\Zn,T)$ is
trivial.
\end{proof}

\begin{theorem}\label{T:Zn-cocycle}
Let $G$ be a connected Lie group with a $\Zn$-module structure,
$\sigma$ a generator of $\Zn$. Let $Z=\{g\in
G|g\sigma(g)\cdots\sigma^{n-1}(g)=e\}$. Then we have\\
(i) $Z$ has finitely many connected components, each of which is a
closed submanifold of $G$;\\
(ii) Under the identification $Z^1(\Zn,G)\cong Z$, each cohomology
class is a connected component of $Z$, that is, $H^1(\Zn,G)$
coincides with the set of all connected components of $Z$.
\end{theorem}

Theorem \ref{T:Zn-cocycle} is equivalent to Theorem
\ref{T:A-cocycle} for the case that $A=\Zn$. In fact, the map
$\rho:G\rightarrow G^n$ defined by
$\rho(g)=(e,g,g\sigma(g),\cdots,g\sigma(g)\cdots\sigma^{n-2}(g))$
is an embedding, and $\rho(Z)=Z^1(\Zn,G)$ under the natural
identification $G^n\cong\Map(\Zn,G)$. But we would like to provide
another proof of Theorem \ref{T:Zn-cocycle}, without the use of
Weil's Lemma.

\begin{proof}[Proof of Theorem \ref{T:Zn-cocycle}]
First we note that the formula
$$\Ad(g\sigma(g)\cdots\sigma^{k-1}(g))d\sigma^k=(\Ad(g)d\sigma)^k$$
holds for every $g\in G$ and every positive integer $k$. In fact,
with respect to the $\ZZ$-module structure on $G$ defined by
$(m,h)\mapsto\sigma^m(h)$, there is a cocycle
$\alpha:\ZZ\rightarrow G$ such that
$\alpha(m)=g\sigma(g)\cdots\sigma^{m-1}(g)$ for $m>0$. Twisting
the $\ZZ$-module $G$ using the cocycle $\alpha$, we know that
$(m,h)\mapsto\alpha(m)\sigma^m(h)\alpha(m)^{-1}$ is an action of
$\ZZ$ on $G$. Differentiating both sides of the equation
$\alpha(k)\sigma^k(h)\alpha(k)^{-1}=(g\sigma(h)g^{-1})^k$, we get
the desired formula.

Now let $z\in Z$. Denote $\ba=\ker(\Ad(z)d\sigma-1)$,
$\ab=\mathrm{Im}(\Ad(z)d\sigma-1)$. By the above formula, we have
$(\Ad(z)d\sigma)^n=1$. So $\Ad(z)d\sigma$ is semisimple, and then
$\g=\ba\oplus\ab$. We construct a smooth map
$F_z:\ba\oplus\ab\rightarrow G$ by
$$F_z(X,Y)=e^{-Y}e^Xz\sigma(e^Y).$$ Then it is easy to compute the
differential $dF_z(0,0):\ba\oplus\ab\rightarrow T_zG$ as
$$dF_z(0,0)(X,Y)=(dr_z)_e(X+(\Ad(z)d\sigma-1)(Y)),$$ where $r_z$
is the right translation on $G$ induced by $z$. Since
$\Ad(z)d\sigma$ is semisimple, the restriction of
$\Ad(z)d\sigma-1$ on $\ab=\mathrm{Im}(\Ad(z)d\sigma-1)$ is a
linear automorphism. Hence $dF_z(0,0)$ is a linear isomorphism,
and then $F_z$ is a local diffeomorphism at $(0,0)$. Choose an
open neighborhood $U_1$ of $0\in\ba$ and an open neighborhood
$U_2$ of $0\in\ab$ such that the restriction of $F_z$ to
$U_1\times U_2\subset\ba\oplus\ab$ is a diffeomorphism onto an
open neighborhood $V=F_z(U_1\times U_2)$ of $z\in G$. Shrinking
$U_1$ if necessary, we may assume that $X\in U_1$, $e^{nX}=e$
implies $X=0$. For $X\in U_1, Y\in U_2$, we have
\begin{align*}
&F_z(X,Y)\sigma(F_z(X,Y))\cdots\sigma^{n-1}(F_z(X,Y))\\
=&e^{-Y}e^Xz\sigma(e^Xz)\cdots\sigma^{n-1}(e^Xz)e^Y\\
=&e^{-Y}e^Xe^{\Ad(z)d\sigma(X)}e^{\Ad(z\sigma(z))d\sigma^2(X)}\cdots
e^{\Ad(z\sigma(z)\cdots\sigma^{n-2}(z))d\sigma^{n-1}(X)}e^Y\\
=&e^{-Y}e^Xe^{\Ad(z)d\sigma(X)}e^{(\Ad(z)d\sigma)^2(X)}\cdots
e^{(\Ad(z)d\sigma)^{n-1}(X)}e^Y\\
=&e^{-Y}e^{nX}e^Y.
\end{align*}
So $F_z(X,Y)\in Z\Leftrightarrow X=0$, that is, $$Z\cap
V=F_z(\{0\}\times U_2)=\{e^{-Y}z\sigma(e^Y)|Y\in U_2\}.$$ This
show that each connected component of $Z$ is a submanifold of $G$,
which is necessary closed by the definition of $Z$, and that each
cohomology class is open in $Z$. So each cohomology class is also
closed in $Z$. But $G$ is connected implies that cohomology
classes are connected. Hence each cohomology class is in fact a
connected component of $Z$. This proves the theorem.
\end{proof}

Hence we get a new proof of the following assertion, which was
first proved in \cite{AW3} using real analytic geometry.

\begin{corollary}[\cite{AW3}, Theorem 1.1]\label{C:E_n}
Let $G$ be a connected Lie group, $n$ a positive integer. Then
each connected component of the set $\{g\in G|g^n=e\}$ is a
conjugacy class of $G$.
\end{corollary}

\begin{proof}
Consider the $\Zn$-module $G$ for which $\Zn$ acts trivially on
$G$ and apply Theorem \ref{T:Zn-cocycle}.
\end{proof}

\begin{remark}
The idea of the construction of the map $F_z$ in the proof of
Theorem \ref{T:Zn-cocycle} is due to Professor K.-H. Neeb.
\end{remark}

We also have

\begin{proposition}\label{P:Zn-coboundary}
Let $G$ be a connected Lie group with a $\Zn$-module structure,
$T$ a maximal compact torus of $G_0^{\Zn}$. Then every connected
component of $H^0(\Zn,G/T)$ is a closed homogeneous submanifold of
$G/T$, and the restriction of the coboundary operator
$\delta:H^0(\Zn,G/T)\rightarrow H^1(\Zn,T)$ to each connected
component of $H^0(\Zn,G/T)$ is constant.
\end{proposition}

\begin{proof}
Similar to the proof of Corollary \ref{C:Z-coboundary}, using the
fact that $H^1(\Zn,T)$ is finite.
\end{proof}

%%%%%%%%%%%%%%%%%%%%%%%%%%%%%%%%%%%%%%%%%%%%%%%%%%%%%%%  Section 6

\section{Twisted conjugate actions of Lie groups}

Let $G$ be a connected Lie group, and let $\sigma$ be an
automorphism of $G$. The twisted conjugate action of $G$ on itself
associated with $\sigma$ is defined by
$\tau_g(h)=gh\sigma(g)^{-1}$. An orbit of the twisted conjugate
action is called a twisted orbit. If $\sigma$ is the identity, the
associated twisted conjugate action is the adjoint action of $G$.

The notion of twisted conjugate action was first introduced in
\cite{AW1,AW2}, where the authors considered the problem of
embedding symmetric spaces into Lie groups. For a connected Lie
group $G$ with an involution $\sigma$ such that $G^\sigma$ is
compact, $G/G^\sigma$ has a structure of Riemannian symmetric
space. It was proved in \cite{AW1} that every connected component
of the set $R=\{g\in G|\sigma(g)=g^{-1}\}$ is a closed submanifold
of $G$. This is in fact a particular case of item (i) of Theorem
\ref{T:Zn-cocycle} in this paper. It was also proved in \cite{AW1}
that the identity component of $R$ coincides with $P=\exp(\ap)$,
where $\ap$ is the eigenspace of $d\sigma$ in the Lie algebra $\g$
of $G$ with eigenvalue $-1$. This gave an isomorphism
$G/G^\sigma\cong P$.

In this section we consider general properties of twisted
conjugate actions of Lie groups. Let $G$ be a Lie group acts
smoothly on two manifolds $M_1$ and $M_2$ by $\rho_i:G\times
M_i\rightarrow M_i$ $(i=1,2)$. The two actions $\rho_1$ and
$\rho_2$ are \emph{equivalent} if there exists an equivariant
diffeomorphism $\varphi:M_1\rightarrow M_2$, that is,
$\varphi(\rho_1(g,x))=\rho_2(g,\varphi(x))$ for every $g\in G$ and
$x\in M_1$.

For a Lie group $G$, we denote the automorphism group and the
inner automorphism group of $G$ by $\Aut(G)$ and $\Inn(G)$,
respectively.

\begin{proposition}\label{P:twisted-iso}
Let $G$ be a connected Lie group. Then automorphisms of $G$ belong
to the same coset of $\Inn(G)$ in $\Aut(G)$ induce equivalent
twisted conjugate actions.
\end{proposition}

\begin{proof}
Let $\sigma, \sigma'\in\Aut(G)$, and denote the twisted conjugate
actions associated with $\sigma$ and $\sigma'$ be $\tau$ and
$\tau'$, respectively. Suppose there exists $g\in G$ such that
$\sigma=\Inn(g)\circ\sigma'$, where $\Inn(g)$ is the inner
automorphism of $G$ induced by $g$. Then
$r_g\circ\tau_h(k)=hk\sigma(h)^{-1}g=hk(g\sigma'(h)g^{-1})^{-1}g
=hkg\sigma'(h)^{-1}=\tau'_h\circ r_g(k)$ for every $h,k\in G$.
Hence the right multiplication $r_g:G\rightarrow G$ is
equivariant. This proves the proposition.
\end{proof}

Proposition \ref{P:twisted-iso} implies that twisted conjugate
actions associated with inner automorphisms are equivalent to the
adjoint action. In particular, if $G$ is semisimple, there are
only finitely many equivalent types of actions appear among
twisted conjugate actions, due to the fact that $\Aut(G)/\Inn(G)$
is finite for $G$ semisimple. Even this, new types of actions do
appear.

\begin{proposition}\label{P:fixedpoint}
Let $G$ be a connected Lie group with an automorphism $\sigma$. If
$\sigma$ is not an inner automorphism, then the twisted conjugate
action $\tau$ associated with $\sigma$ has no fixed point.
\end{proposition}

\begin{proof}
Suppose $\tau$ has a fixed point $h\in G$, that is,
$\tau_g(h)=gh\sigma(g)^{-1}=h$ for every $g\in G$. Then
$\sigma(g)=h^{-1}gh$ for every $g\in G$, a contradiction.
\end{proof}

Since the identity element is a fixed point of the adjoint action,
Proposition \ref{P:fixedpoint} implies that twisted conjugate
actions associated with automorphisms which are not inner are not
equivalent to the adjoint action.

Now we use some results about nonabelian cohomology of cyclic
groups proved in previous sections to prove some properties on the
structure of orbits of twisted conjugate actions.

\begin{theorem}\label{T:compact-twisted}
Let $G$ be a connected compact Lie group, $\sigma\in\Aut(G)$. Let
$T$ be a maximal torus of $G_0^\sigma$. Then the
following statements are equivalent.\\
(i) $\sigma$ is $1$-semisimple;\\
(ii) $O\cap T$ is nonempty for every twisted orbit $O$ of $G$
associated
with $\sigma$;\\
(iii) $O_e\cap T$ is finite, where $O_e$ is the twisted orbit
through the identity element $e$;\\
(iv) $O\cap T$ is finite for every twisted orbit $O$ of $G$
associated with $\sigma$.
\end{theorem}

\begin{proof}
Define the $\ZZ$-module structure on $G$ by
$(m,g)\mapsto\sigma^m(g)$. Then by the relation between twisted
conjugate actions and nonabelian cohomology (see Section 2), (ii)
is equivalent to that $H^1(\ZZ,T)\rightarrow H^1(\ZZ,G)$ is
surjective, (iii) is equivalent to that the kernel of
$H^1(\ZZ,T)\rightarrow H^1(\ZZ,G)$ is finite, and (iv) is
equivalent to that all fibers of $H^1(\ZZ,T)\rightarrow
H^1(\ZZ,G)$ are finite. Hence the theorem follows immediately from
Theorem \ref{T:Z-sur}.
\end{proof}

\begin{theorem}\label{T:noncompact-twisted}
Let $G$ be a connected Lie group, $\sigma$ an automorphism of $G$
of finite order. Let $Z=\{g\in
G|g\sigma(g)\cdots\sigma^{n-1}(g)=e\}$, where $n$ is a positive
integer which is divisible
by the order of $\sigma$. Then\\
(i) $Z$ has finitely many connected components, each of which is a
twisted orbit;\\
(ii) For every maximal compact torus $T$ of $G_0^\sigma$ and every
twisted orbit $O$ contained in $Z$, $O\cap T$ is nonempty and
finite;\\
(iii) There exists a maximal compact subgroup $K$ of $G$ which is
$\sigma$-invariant. For every such $K$ and every twisted orbit $O$
of $G$ contained in $Z$, $O\cap K$ is a twisted orbit of $K$.
\end{theorem}

\begin{proof}
Define the $\Zn$-module structure on $G$ by
$(m,g)\mapsto\sigma^m(g)$. Then under the identification
$Z^1(\Zn,G)\cong Z$ as explained in Section 2, every cohomology
class is a twisted orbit in $Z$ . So (i) follows from Theorem
\ref{T:Zn-cocycle}, (ii) follows from Theorem \ref{T:Zn-sur} and
the fact that $H^1(\Zn,T)$ is finite, and (iii) follows from
Theorem \ref{T:A-bi}.
\end{proof}

\end{document}